\documentstyle[mssymb,11pt,twoside]{article}
\newcommand{\ZZ}{\Bbb Z}

\newcommand{\RR}{\Bbb R}
\newcommand{\CC}{\Bbb C}
\newcommand{\qhyp}[5]{\mbox{}_{#1}\phi_{#2}\left(\left.\begin{array}{c}
{#3}\\{#4}\end{array}\right| q ; {#5}\right)}
\newcommand{\qqhyp}[5]{\mbox{}_{#1}\phi_{#2}\left(\left.\begin{array}{c}
{#3}\\{#4}\end{array}\right| q^2 ; {#5}\right)}

\begin{document}

\title{On the zeros of the Hahn-Exton $q$-Bessel function and 
associated $q$-Lommel polynomials}
\author{H.T. Koelink\thanks{Van Swietenstraat 3, 2334 EA Leiden,
The
Netherlands} \and R.F. Swarttouw\thanks{Vrije Universiteit,
Faculteit
Wiskunde en Informatica, De Boelelaan 1081, 1081 HV  Amsterdam, The
Netherlands}}
\date{\ }

\maketitle

\section{Introduction}

\markboth{\hfill{\sc Preprint}\hfill}{\hfill{\sc Preprint}\hfill}
\pagestyle{myheadings}

\def\theequation{1.\arabic{equation}}
\setcounter{equation}{0}

For the Bessel function
\begin{equation}
\label{bessel} J_{\nu}(z) = \sum\limits_{k=0}^{\infty}
\frac{(-1)^k \left( \frac{z}{2} \right)^{\nu+2k}}{k!
\Gamma(\nu+1+k)}
\end{equation}
there exist several $q$-analogues. The oldest $q$-analogues of the
Bessel function were introduced by F.H. Jackson at the beginning 
of this century, see M.E.H. Ismail \cite{Is1} for the appropriate 
references.
Another $q$-analogue of the Bessel function has been introduced by 
W. Hahn in a special case and by H. Exton in full generality, see
R.F. Swarttouw \cite{Sw1} for a historic overview.

Here we concentrate on properties of the Hahn-Exton $q$-Bessel
function
and in particular on its zeros and the associated $q$-Lommel 
polynomials. Ismail has proved very satisfactory results on this
subject
for the Jackson $q$-Bessel functions, cf. \cite{Is1}. In
particular, he
proves orthogonality relations for the associated $q$-Lommel
polynomials.

In section 2 we present the definition of the Hahn-Exton $q$-Bessel
function and some of its properties. The zeros of the Hahn-Exton
$q$-Bessel function of order $\nu>-1$ are the subject of section 3.
The
proofs of the statements in this section rest on the evaluation of
a
$q$-integral closely related to the Fourier-Bessel orthogonality
relations for the Hahn-Exton $q$-Bessel function. The results of
section 3 are in accordance with the results on the zeros of the
Bessel
function. Section 4 deals with
two kinds of associated $q$-Lommel polynomials for which we present
explicit forms as well as some properties. However, the first type
does not give rise to orthogonal polynomials, while the second type
provides polynomials, which are closely related to the orthogonal
modified $q$-Lommel polynomials found by M.E.H. Ismail.

Comparison of our results with the results for the Jackson
$q$-Bessel
functions shows that the Hahn-Exton $q$-Bessel function is less
similar
to the Bessel function in this respect than the Jackson $q$-Bessel
function, since we are not able to prove new orthogonality
relations for
the $q$-Lommel polynomials associated with the Hahn-Exton
$q$-Bessel
function. This seems in support of
M. Rahman \cite{Ra1}, who favours the Jackson $q$-Bessel function
over
the Hahn-Exton $q$-Bessel function. However, it should be noted
that
the proofs concerning the zeros of the Hahn-Exton $q$-Bessel
function
are much simpler and more direct than those for the Jackson
$q$-Bessel
function.

In our opinion both $q$-analogues of the Bessel function are
interesting functions and possess nice properties. Moreover, from
the
harmonic analysis on the quantum group of plane motions, cf.
H.T. Koelink [5, \S6.7], it follows that there exists a common
generalisation of the Jackson and Hahn-Exton $q$-Bessel function.

\section{The Hahn-Exton $q$-Bessel function}

\markboth{\hfill{\sc Preprint}\hfill}{\hfill{\sc Preprint}\hfill}
\pagestyle{myheadings}

\def\theequation{2.\arabic{equation}}
\setcounter{equation}{0}

In this section we present the Hahn-Exton $q$-Bessel function and
some
of its properties. References for these results are T.H.
Koornwinder and
R.F. Swarttouw \cite{Ko1} and R.F. Swarttouw \cite{Sw1}. The
$q$-Bessel
functions
are defined in terms of basic hypergeometric series, which we will
briefly recall. More information on basic hypergeometric series can
be
found in the book by G. Gasper and M. Rahman [3, Chapter 1].

\vspace{5mm}

We fix $0 < q < 1$ and we define the $q$-shifted factorials
\begin{eqnarray}
\label{qshift} (a;q)_k = \left\{ \begin{array}{ll}
1 &\mbox{ if } k=0\\[1ex]
(1-a)(1-aq)\cdots(1-aq^{k-1}) &\mbox{ if } k \geq 1
\end{array}\right.
\end{eqnarray}
for arbitrary $a\in\CC$ and $k\in\ZZ_+$. The limit
$k\rightarrow\infty$
is well defined and yields
\begin{eqnarray}
\label{qshiftinf}(a;q)_{\infty} = \lim\limits_{k\rightarrow\infty}
(a;q)_k.
\end{eqnarray}
The $q$-hypergeometric series (or basic hypergeometric) series
$_r\phi_s$ is defined by
\begin{eqnarray}
\qhyp{r}{s}{a_1, \dots, a_r}{b_1, \dots, b_s}{z} =
\sum\limits_{k=0}
^{\infty}\frac{(a_1;q)_k(a_2;q)_k\cdots(a_r;q)_k}{(q;q)_k(b_1;q)_
k\cdots
(b_s;q)_k} \left( (-1)^kq^{\frac{1}{2}k(k-1)}\right)^{1+s-r} z^k
\label{qhyper}
\end{eqnarray}
whenever the series converges.

\vspace{5mm}

The Hahn-Exton $q$-Bessel function $J_{\nu}(x;q)$ of order $\nu$ is
defined as
\begin{eqnarray}
J_{\nu}(x;q) = x^{\nu}
\frac{(q^{\nu+1};q)_{\infty}}{(q;q)_{\infty}}
\qhyp{1}{1}{0}{q^{\nu+1}}{qx^2}.\label{heq}
\end{eqnarray}
From (\ref{qhyper}) we see that $x^{-\nu} J_{\nu}(x;q^2)$
defines a non-zero analytic function on $\CC$, so that
$J_{\nu}(x;q)$
is analytic in $\CC\backslash\{0\}$.
The Hahn-Exton $q$-Bessel function is a $q$-analogue of the Bessel
function, since
$J_{\nu}((1-q)x;q)$ tends to the Bessel function
$J_{\nu}(2x)$ of order $\nu$ as $q\uparrow 1$, cf. [9, Appendix A].
The
standard text on Bessel functions is the classic by G.N. Watson
\cite{Wa1}.

\vspace{5mm}

Next we introduce two concepts of $q$-analysis: the $q$-derivative
and
the $q$-integral, cf. [3, Chapter 1]. The $q$-derivative of a
function
$f$ is defined by
\begin{eqnarray}
\left( D_q f\right)(x) = \frac{f(x) - f(qx)}{(1-q)x},\hspace{5mm}
x \neq 0.\label{qaf}
\end{eqnarray}
Note that the $\left( D_q f \right)(x)$ tends to $f'(x)$ for
$q\uparrow
1$ whenever $f$ is differentiable at $x$. The product rule for the
$q$-derivative is
\begin{eqnarray}
\label{qprod}\left(D_q fg\right)(x) = f(x) \left(D_q g\right)(x) +
g(qx) \left(D_q f\right)(x).
\end{eqnarray}
The $q$-integral of a function $f$ is defined by
\begin{eqnarray}
\label{qint}\int\limits_0^z f(x) d_qx = (1-q)z \sum\limits_{k=0}^
{\infty} f(zq^k) q^k, \hspace{5mm} z > 0.
\end{eqnarray}
For a Riemann-integrable function $f$ this expression tends to
$\int_0^z f(x)\, dx$ as $q\uparrow 1$. It is easily checked that
for
continuous $f$
\begin{eqnarray}
\int_0^z \left(D_q f\right)(x)\, d_qx = f(z) - f(0).\label{intdif}
\end{eqnarray}
The $q$-product rule leads to the following $q$-partial integration
rule
\begin{eqnarray}
\int\limits_0^z f(qx) \left( D_q g\right)(x) \, d_q x = f(z)g(z) -
f(0)g(0)
-\int\limits_0^z \left( D_q f \right)(x) g(x) \, d_q
x.\label{qpart}
\end{eqnarray}

\vspace{5mm}

We will use the following $q$-derivatives, cf. [9, (3.2.22),
(3.2.19)]:
\begin{eqnarray}
& &D_q \left[ (\cdot)^{\nu} J_{\nu}(\cdot \, ;q^2)\right](x) =
\frac{x^{\nu}}{1-q} J_{\nu-1}(x;q^2),\label{af1}\\[1ex]
& &D_q \left[ (\cdot)^{-\nu} J_{\nu}(\cdot \, ;q^2)\right](x) =
- \frac{q^{1-\nu} x^{-\nu}}{1-q} J_{\nu+1}(xq;q^2).\label{af2}
\end{eqnarray}
These formulas are $q$-analogues of the case $m=1$ of
[11, 3.2(5), 3.2(6)].

\vspace{5mm}

The second order differential equation for the Bessel function, cf.
[11, 3.2(1)], has the following second order $q$-difference
equation as
a $q$-analogue for the Hahn-Exton $q$-Bessel function, cf. [9,
(4.3.1)]:
\begin{eqnarray}
\label{qdifeq}J_{\nu}(xq^2;q^2) + q^{-\nu} \left( x^2q^2 - 1 -
q^{2\nu}
\right) J_{\nu}(xq;q^2) + J_{\nu}(x;q^2) = 0.
\end{eqnarray}
It will also be useful to have the following relations at hand.
They
involve Hahn-Exton $q$-Bessel functions of several orders and can
be
found in [9, (3.2.15), (3.2.18), (3.2.21)]:
\begin{eqnarray}
& &J_{\nu+1}(x;q^2) = \left( \frac{1-q^{2\nu}}{x} + x \right)
J_{\nu}
(x;q^2) - J_{\nu-1}(x;q^2),\label{mix1}\\[1ex]
& &J_{\nu+1}(xq;q^2) = q^{-\nu-1} \left( \frac{1-q^{2\nu}}{x}
J_{\nu}
(x;q^2) - J_{\nu-1}(x;q^2) \right).\label{mix2}
\end{eqnarray}

\section{On the zeros of the Hahn-Exton $q$-Bessel function}
\def\theequation{3.\arabic{equation}}
\setcounter{equation}{0}

In this section we prove some results on the zeros of the
Hahn-Exton
$q$-Bessel function. The proofs rest on an explicit evaluation of
a
$q$-integral related to the Fourier-Bessel orthogonality relations
for
the Hahn-Exton $q$-Bessel function and on some formulas relating 
Hahn-Exton $q$-Bessel functions of different order. At the end of
this
section we also present a result on the zeros of the $q$-derivative
of
the Hahn-Exton $q$-Bessel function.

\vspace{5mm}

The starting point in the derivation of the results on the zeros of
the
Hahn-Exton $q$-Bessel function is the following proposition. It is
closely related to the Fourier-Bessel orthogonality relations for
the
Hahn-Exton $q$-Bessel functions, originally proved by H. Exton (see
\cite{Ex1} and \cite{Ex2}) using Sturmian methods.

\vspace{5mm}

{\bf Proposition 3.1.} {\sl For $\Re e(\nu) > -1, \, z > 0$ and 
$a,b\in\CC
\backslash\{0\}$ we have}
\begin{eqnarray*}
& &\left(a^2 - b^2 \right) \int\limits_0^z x J_{\nu}(aqx;q^2)
J_{\nu}
(bqx;q^2) \, d_q x\\[1ex]
& &\hspace{1cm}= (1-q) q^{\nu - 1} z \left( a J_{\nu+1}(aqz;q^2)
J_{\nu}(bz;q^2) - b J_{\nu+1}(bqz;q^2) J_{\nu}(az;q^2) \right).
\end{eqnarray*}

\vspace{5mm}

{\bf Proof:} Use the $q$-partial integration rule (\ref{qpart}) and
the
formulas (\ref{af1}) and (\ref{af2}) to obtain
\begin{eqnarray}
& &\int\limits_0^z x J_{\nu}(aqx;q^2) J_{\nu}
(bqx;q^2) \, d_qx\label{proof1}\\[1ex]
\hspace{1cm}& &= \frac{(1-q)}{b} q^{\nu - 1} z  J_{\nu}(az;q^2)
J_{\nu+1}(bqz;q^2) + \frac{a}{b} \int\limits_0^z x
J_{\nu+1}(aqx;q^2)
J_{\nu+1}(bqx;q^2)\, d_qx.\nonumber
\end{eqnarray}
The restriction $\Re e(\nu) > -1$ is necessary to ensure the
absolute
convergence
of the series defining the $q$-integral, cf. (\ref{qint}), since
the
Hahn-Exton $q$-Bessel function $J_{\nu}(x;q^2)$ behaves like
$x^{\nu}$
times a constant for small $x$. Interchanging
$a$ and $b$ in (\ref{proof1}) yields a set of two equations, which
can
be solved easily. $\Box$

\vspace{2mm}

The proof of Proposition 3.1 is also well known for the Bessel
function, see e.g. I.N. Sneddon [8, Chapter 2].

\vspace{5mm}

{\bf Corollary 3.2.} {\sl The zeros of $J_{\nu}(\cdot\, ;q^2)$,
with
$\nu >
-1$, are real.}

\vspace{5mm}

{\bf Proof:} Suppose $a\neq 0$ is a zero of $J_{\nu}(\cdot\,
;q^2)$.
Since $\nu$ is real we have
\begin{eqnarray*}
J_{\nu}\left(\bar{a};q^2\right) = \overline{J_{\nu}(a;q^2)} = 0.
\end{eqnarray*}
Proposition 3.1 with $z = 1$ and $b = \bar{a}$ yields
\begin{eqnarray}
\left( a^2 - \bar{a}^2 \right) \int\limits_0^1 x \left|
J_{\nu}(aqx;q^2)
\right|^2 \, d_qx = 0.\label{abar}
\end{eqnarray}
Now $a^2 = \bar{a}^2$ if and only if $a\in\RR$ or $a\in i\RR$, so
that
in all other cases the $q$-integral in (\ref{abar}) is zero. Using
the
definition of the $q$-integral (\ref{qint}) we get
\begin{eqnarray*}
J_{\nu}\left(aq^{k+1};q^2 \right) = 0, \hspace{2mm} k\in\ZZ_+,
\end{eqnarray*}
and this implies that $x^{-\nu} J_{\nu}(x;q^2)$
is identically zero,
since it defines an analytic function on $\CC$. Hence,
$J_{\nu}(\cdot\,
;q^2) \equiv 0$.

Finally, we have to show that $ia$, with $a\in\RR$, is not a zero
of
$J_{\nu}(\cdot\, ;q^2)$. Now by (\ref{qhyper}) and (\ref{heq})
\begin{eqnarray*}
J_{\nu}(ia;q^2) = (ia)^{\nu}
\frac{(q^{2\nu+2};q^2)_{\infty}}{(q^2;q^2)
_{\infty}} \sum\limits_{k=0}^{\infty} \frac{q^{k(k+1)} a^{2k}}
{(q^2;q^2)_k
(q^{2\nu+2};q^2)_k}
\end{eqnarray*}
and for $\nu > -1$ this expression cannot be zero, since every term
in
the sum is positive. $\Box$

\vspace{5mm}

From the series representation (\ref{heq}) for the Hahn-Exton
$q$-Bessel
function it follows that if $a$ is a zero of $J_{\nu}(\cdot\,
;q^2)$,
then
$-a$ is also a zero of $J_{\nu}(\cdot\, ;q^2)$. Hence we will
restrict
ourselves to the positive zeros of the Hahn-Exton $q$-Bessel
function
of order $\nu > -1$.

\vspace{5mm}

To obtain an expression for the $q$-integral in Proposition 3.1
with
$a=b$, we use l'H\^{o}pital's rule. The result is
\begin{eqnarray*}
\int\limits_0^z x \left( J_{\nu}(axq;q^2)\right)^2 \, d_qx &=&
\frac
{(1-q)q^{\nu-1}z}{-2a}\left( az J_{\nu+1}(aqz;q^2) J_{\nu}'(az;q^2)
\right.\\[1ex]
& &\left. - J_{\nu+1}(aqz;q^2)J_{\nu}(az;q^2) -aqz
J_{\nu+1}'(aqz;q^2)
J_{\nu}(az;q^2)\right).
\end{eqnarray*}
This formula simplifies to
\begin{eqnarray}
\int\limits_0^1 x \left(J_{\nu}(axq;q^2) \right)^2 \, d_qx =
-\frac{1}
{2}(1-q)q^{\nu-1} J_{\nu+1}(aq;q^2)J_{\nu}'(a;q^2)\label{simp}
\end{eqnarray}
for $z=1$ and $a\neq 0$ a (real) zero of $J_{\nu}(\cdot\,;q^2)$.

\vspace{5mm}

{\bf Lemma 3.3.} {\sl The non-zero (real) zeros of
$J_{\nu}(\cdot\,;q^2)$,
with $\nu > -1$, are simple zeros.}

\vspace{5mm}

{\bf Proof:} Let $a$ be a non-zero (real) zero of
$J_{\nu}(\cdot\, ;q^2)$,
with $\nu > -1$. The integral
\begin{eqnarray*}
\int\limits_0^1 x \left| J_{\nu}(axq;q^2) \right|^2 \, d_qx =
\int\limits
_0^1 x \left( J_{\nu}(axq;q^2) \right)^2 d_qx
\end{eqnarray*}
is strictly positive. (If it were zero, this would imply that the
Hahn-Exton $q$-Bessel function is identically zero as in the proof
of
Corollary 3.2.) Hence, (\ref{simp}) implies that $J_{\nu}'(a;q^2)
\neq 0$, which proves the lemma. $\Box$

\vspace{2mm}

In the proof of this lemma we explicitly use the (usual) derivative
of
the Hahn-Exton $q$-Bessel function in (\ref{simp}). It shows that
the
intermingling of (ordinary) analysis and $q$-analysis may be
fruitful
and should not be opposed to as M. Rahman \cite{Ra1} proposes.

\vspace{5mm}

We now come to one of the main results of this section, which shows
that
the zeros of the Hahn-Exton $q$-Bessel function behave like the
zeros
of the classical Bessel function, cf. Watson [11, \S15.3].

\vspace{5mm}

{\bf Theorem 3.4.} {\sl The Hahn-Exton $q$-Bessel function of order
$\nu > -1$ has a countably infinite number of positive simple
zeros.}

\vspace{5mm}

{\bf Proof:} In view of Corollary 3.2 and Lemma 3.3 it remains to
show
that $J_{\nu}(\cdot\,;q^2)$, with $\nu > -1$, has a countably
infinite
number of zeros. We prove that the assumption
$J_{\nu}(\cdot\,;q^2)$,
with $\nu > -1$, has a finite number of real positive zeros, leads
to
a contradiction.

Assume that $a > 0$ is the largest positive zero of
$J_{\nu}(\cdot\,;q^2)$
and choose $x > 0$ satisfying
\begin{eqnarray*}
x^2q^2 - 1 -q^{2\nu} > 0,\hspace{5mm}\mbox{ and } \hspace{5mm} q^2x
> a.
\end{eqnarray*}
Hence, $J_{\nu}(xq^2;q^2), J_{\nu}(xq;q^2)$ and $J_{\nu}(x;q^2)$
are
non-zero real numbers of the same sign and
\begin{eqnarray*}
J_{\nu}(xq^2;q^2) + q^{-\nu} \left( q^2x^2 - 1 -q^{2\nu} \right)
J_{\nu}
(xq;q^2) + J_{\nu}(x;q^2)
\end{eqnarray*}
is a non-zero real number. This is in contradiction with the second
order $q$-difference equation (\ref{qdifeq}) for the Hahn-Exton
$q$-Bessel function. $\Box$

\vspace{5mm}

Theorem 3.4 shows that we can order the positive zeros of $J_{\nu}
(\cdot\,;q^2)$, with $\nu > -1$, as
\begin{eqnarray}
\label{volg}0 < j_1^{\nu}(q^2) < j_2^{\nu}(q^2) < j_3^{\nu}(q^2) <
\dots .
\end{eqnarray}
The positive zeros of $J_{\nu}(\cdot\, ;q^2)$ will also be denoted
by
$j_n$ or $j_n^{\nu}$. Proposition 3.1 and
relation (\ref{simp}) can be combined to state the
Fourier-Bessel orthogonality relations for the Hahn-Exton
$q$-Bessel
function, cf. \cite{Ex1} and \cite{Ex2}.

\vspace{5mm}

{\bf Proposition 3.5.} {\sl Let $\nu > -1$ and $0 < j_1 < j_2 <
\dots$
be the
positive zeros of the Hahn-Exton $q$-Bessel function
$J_{\nu}(\cdot\, ;
q^2)$, then}
\begin{eqnarray*}
& &\int\limits_0^1 x J_{\nu}(qj_nx;q^2)J_{\nu}(qj_mx;q^2)\, d_q
x\\[1ex]
& &\hspace{1cm}= -\frac{1}{2}(1-q)q^{\nu-1}
J_{\nu+1}(qj_n;q^2) J_{\nu}'(j_n;q^2) \delta_{n,m}\\[1ex]
& &\hspace{1cm}= \frac{1}{2}(1-q)^2q^{\nu-2}\left(D_q J_{\nu}
(\cdot\, ;q^2)\right)(j_n) \, J_{\nu}'(j_n;q^2) \delta_{n,m}\\[1ex]
& &\hspace{1cm}= -\frac{1}{2}(1-q)q^{\nu-2}
j_n^{-1} J_{\nu}(qj_n;q^2) J_{\nu}'(j_n;q^2) \delta_{n,m}\\[1ex]
& &\hspace{1cm}= -\frac{1}{2}(1-q)q^{-2}J_{\nu+1}(j_n;q^2) J_{\nu}'
(j_n;q^2) \delta_{n,m}.
\end{eqnarray*}

\vspace{5mm}

{\bf Proof:} It remains to prove the last three equalities. The
first
one is a consequence of
\begin{eqnarray}
\label{af3}\left[D_q J_{\nu}(\cdot\, ;q^2)\right](x) =
\frac{-q}{1-q}
J_{\nu+1}(xq;q^2) + \frac{1-q^{\nu}}{x(1-q)} J_{\nu}(x;q^2),
\end{eqnarray}
which follows from (\ref{qprod}) and (\ref{af2}).
The second equality follows from the definition of the
$q$-derivative
$D_q$, cf. (\ref{qaf}). The third equality follows from $J_{\nu+1}
(j_n;q^2) = q^{\nu+1} J_{\nu+1}(qj_n;q^2)$, which is a consequence
of
(\ref{mix1}) and (\ref{mix2}). $\Box$

\vspace{5mm}

Next we will derive that the zeros of the Hahn-Exton $q$-Bessel
functions of order $\nu$ and $\nu+1$
are interlaced similarly to the interlacing property of the
zeros of the Bessel functions of order $\nu$ and $\nu+1$, cf.
Watson
[11, \S15.22]. We start with the following lemma.

\vspace{5mm}

{\bf Lemma 3.6.} {\sl Let $\nu > -1$, then $J_{\nu-1}(\cdot\,
;q^2)$ and
$J_{\nu+1}(\cdot\, ;q^2)$ each have at least one zero between two
consecutive positive zeros of $J_{\nu}(\cdot\, ;q^2)$.}

\vspace{5mm}

{\bf Proof:} Consider the function
\begin{eqnarray*}
g(x) = \frac{x^{\nu}}{1-q} J_{\nu-1}(x;q^2) = D_q \left(
(\cdot)^{\nu}
J_{\nu}(\cdot\, ;q^2) \right)(x),
\end{eqnarray*}
where we used (\ref{af1}) as well. If $a$ is a positive zero of
$J_{\nu}(\cdot\, ;q^2)$, then
\begin{eqnarray*}
g(a) = - \frac{(aq)^{\nu} J_{\nu}(aq;q^2)}{(1-q)a}
\end{eqnarray*}
by (\ref{qaf}). It follows from Proposition 3.5 and Lemma 3.3 that
$g(a)
g(b) < 0$ for two consecutive zeros $0 < a < b$ of $J_{\nu}(\cdot\,
;q^2)$. This proves the lemma for $J_{\nu-1}(\cdot\, ;q^2)$.

The other case is reduced to the result for
$J_{\nu-1}(\cdot\,;q^2)$,
since $J_{\nu-1}(a;q^2) = - J_{\nu+1}(a;q^2)$ for a non-zero zero
of
$J_{\nu}(\cdot\,;q^2)$ by (\ref{mix1}). $\Box$

\vspace{5mm}

{\bf Theorem 3.7.} {\sl The positive real zeros of $J_{\nu}(\cdot\,
;q^2)$
and $J_{\nu+1}(\cdot\,;q^2)$, with $\nu > -1$, interlace.
Explicitly,}
\begin{eqnarray*}
0 < j_1^{\nu} < j_1^{\nu+1} < j_2^{\nu} < j_2^{\nu+1} < j_3^{\nu}
<
j_3^{\nu+1} < \dots.
\end{eqnarray*}

\vspace{5mm}

{\bf Proof:} The interlacing of the zeros of $J_{\nu}(\cdot\,;q^2)$
and
$J_{\nu+1}(\cdot\, ;q^2)$ is a straightforward application of Lemma
3.6.
It remains to prove that $j_1^{\nu} < j_1^{\nu+1}$.

First note that the $_1\phi_1$-series in (\ref{heq}) yields $1$ for
$x=0$. Hence, by continuity,
\begin{eqnarray}
\label{jpos}J_{\nu}(x;q^2) > 0, \hspace{5mm} x\in\left( 0,
j_1^{\nu}
\right), \, \nu > -1.
\end{eqnarray}
Thus $J_{\nu}'(j_1^{\nu};q^2) < 0$ and by Lemma 3.3 and Proposition
3.5
we get $J_{\nu+1}(j_1^{\nu};q^2) > 0$. From (\ref{jpos}) it follows
that
there exists an even number of positive zeros of $J_{\nu+1}(\cdot\,
;q^2)$ in the open interval $(0, j_1^{\nu})$. Lemma 3.6 and
$j_1^{\nu}$
being the smallest positive zero of $J_{\nu}(\cdot\, ;q^2)$ imply
that
this number is zero. $\Box$

\vspace{5mm}

An upper bound for the first zero $j^{\nu}_1$ may be obtained as
follows. Take $x=j^{\nu}_1$ in (\ref{qdifeq}) to obtain
\begin{eqnarray*}
J_{\nu}(j^{\nu}_1 q^2;q^2) + q^{-\nu} \left( \left(j^{\nu}_1
\right)^2
q^2 - 1 - q^{2\nu} \right) J_{\nu}(j^{\nu}_1 q;q^2) = 0.
\end{eqnarray*}
Inequality (\ref{jpos}) implies that both $q$-Bessel functions are
positive, so that
\begin{eqnarray*}
\left(j^{\nu}_1 \right)^2 q^2 - 1 - q^{2\nu} < 0
\Longleftrightarrow
j^{\nu}_1 < q^{-1} \sqrt{1+q^{2\nu}}.
\end{eqnarray*}

\vspace{5mm}

Theorem 3.7 shows that for $\nu>-1$ the Hahn-Exton $q$-Bessel
functions
of order $\nu$ and $\nu+1$ do not have zeros in common except
possibly
$0$. This also holds for general order.

\vspace{2mm}

{\bf Proposition 3.8.} {\sl The Hahn-Exton $q$-Bessel functions
$J_{\nu}(\cdot \, ;q^2)$ and $J_{\nu+1}(\cdot \, ;q^2)$ have no
common
zeros except possibly $0$.}

\vspace{2mm}

{\bf Proof:} We argue by contradiction, so let $x\in\CC$ be a
non-zero
zero of both $J_{\nu}(\cdot \, ;q^2)$ and $J_{\nu-1}(\cdot \,
;q^2)$.
Eliminating $J_{\nu-1}(x;q^2)$ from (\ref{mix1}) and (\ref{mix2})
leads
to
\begin{eqnarray*}
J_{\nu+1}(x;q^2) - q^{\nu+1} J_{\nu+1}(xq;q^2) = x J_{\nu}(x;q^2).
\end{eqnarray*}
Hence, $J_{\nu+1}(xq;q^2) = 0$ as well. Next consider the
$q$-difference
equation (\ref{qdifeq}) with $\nu$ replaced by $\nu+1$. Since $x$
and
$xq$ are zeros of $J_{\nu+1}(\cdot \, ;q^2)$ we find
$J_{\nu+1}(xq^2;q^2)=0$. Repeated application of (\ref{qdifeq})
leads
to $J_{\nu+1}(xq^k;q^2)=0$ for all $k\in\ZZ_+$. This implies
$J_{\nu+1}(\cdot \, ;q^2)\equiv 0$, which is the desired
contradiction.
$\Box$

\vspace{5mm}

Finally we state a result on the zeros of the $q$-derivative of the
Hahn-Exton $q$-Bessel function. Specialise $z=1$ in Proposition 3.1
and use (\ref{af3}) to obtain
\begin{eqnarray}
& &\label{spec}(a^2 - b^2) \int\limits_0^1 x J_{\nu}(aqx;q^2)
J_{\nu}
(bqx;q^2) \, d_q x \\[1ex]
& &\hspace{1cm}= (1-q)^2 q^{\nu-2} \left( b \left(D_q
J_{\nu}(\cdot\,;q^2)
\right)(b) J_{\nu}(a;q^2) - a \left(D_q
J_{\nu}(\cdot\,;q^2)\right)(a)
J_{\nu}(b;q^2) \right)\nonumber
\end{eqnarray}
for $\Re e(\nu) > -1$ and
$a,b\in\CC\backslash\{0\}$.

\vspace{5mm}

{\bf Proposition 3.9.} {\sl The non-zero zeros of $D_q
J_{\nu}(\cdot\,
;q^2)$ for $\nu > 0$ are real and simple.}

\vspace{5mm}

{\bf Proof:} In order to prove that the zeros of
$D_q J_{\nu}(\cdot\,;q^2)$ are
real, we proceed as in the proof of Corollary 3.2. In this case
purely imaginary zeros can be excluded for $\nu > 0$.

To see that the zeros are simple, we apply l'H\^{o}pital's rule to
(\ref{spec}) and then take $a$ a non-zero real zero of $D_q J_{\nu}
(\cdot\,;q^2)$. The result is
\begin{eqnarray*}
\int\limits_0^1 x \left| J_{\nu}(aqx;q^2) \right|^2 \, d_q x = -
\frac
{1}{2}(1-q)^2q^{\nu-2} J_{\nu}(a;q^2) \left( D_q J_{\nu}(\cdot\,
;q^2)
\right)'(a).
\end{eqnarray*}
Since the left hand side is non-zero, the result follows.
$\Box$

\section{\bf On the associated $q$-Lommel polynomials.}
\def\theequation{4.\arabic{equation}}
\setcounter{equation}{0}

The Bessel function $J_{\nu}(z)$ satisfies the three term
recurrence
relation
\begin{eqnarray}
\label{recurbes}J_{\nu+1}(z) = \frac{2\nu}{z} J_{\nu}(z) -
J_{\nu-1}(z).
\end{eqnarray}
If we iterate relation (\ref{recurbes}), we can express
$J_{\nu+m}(z)$,
with $m\in\ZZ_+$, in terms of $J_{\nu}(z)$ and $J_{\nu-1}(z)$, with
coefficients that are polynomials in $\frac{1}{z}$. Indeed we have
\begin{equation}
\label{recurmbes}J_{\nu+m}(z) = R_{m,\nu}(z)J_{\nu}(z) -
R_{m-1,\nu+1}
(z) J_{\nu-1}(z),
\end{equation}
where $R_{m,\nu}(z)$ are the Lommel polynomials, see Watson [11,
\S9.6].
These polynomials satisfy the three term recurrence relation, cf.
[11, 9.63(2)]
\begin{equation}
\label{recurlom}R_{m+1,\nu}(z) = \frac{2(\nu+m)}{z} R_{m,\nu}(z) -
R_{m-1,\nu}(z),
\end{equation}
with $R_{0,\nu}(z) = 1$ and $R_{1,\nu}(z) = \frac{2\nu}{z}$.
Usually
a related set of polynomials, the so-called modified Lommel
polynomials, is defined by
\begin{eqnarray*}
h_{m,\nu}(z) = R_{m,\nu}\left(\frac{1}{z}\right).
\end{eqnarray*}
It is clear that $h_{m,\nu}(z)$ is a polynomial in $z$ of degree
$m$,
which satisfies the three term recurrence relation
\begin{eqnarray*}
\label{recurmod}h_{m+1,\nu}(z) = 2z(\nu+m) h_{m,\nu}(z) -
h_{m-1,\nu}(z).
\end{eqnarray*}
By Favard's theorem, the modified Lommel polynomials are orthogonal
with respect to a positive measure. The orthogonality measure is a
discrete measure with weights at $\frac{1}{j_n^{\nu}}$, where
$j_n^{\nu}$ are the zeros of the Bessel function $J_{\nu}(z)$.

\vspace{5mm}

The asymptotic behaviour of the Lommel polynomials $R_{m,\nu}$ is
related to the Bessel function by Hurwitz's formula, cf. [11,
9.65(1)]
\begin{eqnarray}
\label{hur} \frac{\left(\frac{1}{2}z\right)^{\nu+m} R_{m,\nu+1}(z)}
{\Gamma(\nu+m+1)} \stackrel{m\rightarrow\infty}{\longrightarrow}
J_{\nu}(z).
\end{eqnarray}

In the next subsections we will derive and investigate two
different
$q$-analogues
of the Lommel polynomials. The first type follows from the
recurrence
relation (\ref{mix1}). The second type is a result of the
difference-recurrence relation (\ref{mix2}).

\subsection{$q$-Lommel polynomials associated with (2.13)}

In this subsection we determine $q$-analogues of the Lommel
polynomials
arising from the recurrence relation (\ref{mix1}). We will present
an
explicit formula, a generating function and an analogue of
Hurwitz's
formula (\ref{hur}).
In order to derive $q$-analogues of the results mentioned in the
previous section, we introduce a second Hahn-Exton $q$-Bessel
function.
Define
\begin{eqnarray}
{\cal J}_{\nu}(x;q) &=& e^{i\nu\pi} q^{-\frac{1}{2}\nu}J_{-\nu}
(xq^{-\frac{1}{2}\nu};q)\label{calj},\\[1ex]
&=& e^{i\nu\pi} \frac{(q^{-\nu +1};q)_{\infty}
}{(q;q)_{\infty}}x^{-\nu}q^{\frac{1}{2}\nu
(\nu-1)}\sum\limits_{k=0}
^{\infty}\frac{(-1)^kq^{\frac{1}{2}k(k+1)}x^{2k}q^{-\nu k}}
{(q^{-\nu +1};q)_k(q;q)_k}\nonumber
\end{eqnarray}
Relations (\ref{af1}) and (\ref{af2}) give us, also in combination
with the definition (\ref{calj}) above, the next formulas:
\begin{eqnarray}
\label{drj1}& &J_{\nu}(xq^{\frac{1}{2}};q) = q^{\frac{1}{2}\nu}
J_{\nu}(x;q) + x q^{\frac{1}{2}}J_{\nu
+1}(xq^{\frac{1}{2}};q),\\[1ex]
\label{drjj1}& &{\cal J}_{\nu}(xq^{\frac{1}{2}};q) = q^{\frac{1}{2}
\nu}{\cal J}_{\nu}(x;q) + x q^{\frac{1}{2}}{\cal J}_{\nu
+1}(xq^{\frac
{1}{2}};q),\\[1ex]
\label{drj2}& &J_{\nu}(xq^{\frac{1}{2}};q) =
q^{-\frac{1}{2}\nu}J_{\nu}
(x;q) - xq^{-\frac{1}{2}\nu}J_{\nu -1}(x;q),\\[1ex]
\label{drjj2}& &{\cal J}_{\nu}(xq^{\frac{1}{2}};q) =
q^{-\frac{1}{2}
\nu}{\cal J}_{\nu}(x;q) - xq^{-\frac{1}{2}\nu}{\cal J}_{\nu
-1}(x;q).
\end{eqnarray}

\vspace{5mm}

\noindent Furthermore, ${\cal J}_{\nu}(x;q^2)$ also satisfies the
relations
(\ref{qdifeq}) and (\ref{mix1}), i.e.
\begin{eqnarray}
\label{recjj1}& &\left\{ x + \frac{1-q^{\nu}}{x} \right\}
{\cal J}_{\nu}(x;q) =
{\cal J}_{\nu -1}(x;q) + {\cal J}_{\nu +1}(x;q),\\[1ex]
\label{difj2}& &{\cal J}_{\nu}(xq;q) + q^{-\frac{1}{2}\nu}\left(
x^2q-1-
q^{\nu}\right) {\cal J}_{\nu}(xq^{\frac{1}{2}};q) + {\cal
J}_{\nu}(x;q)
=0.
\end{eqnarray}

\noindent This follows by combining relations (\ref{drjj1}) and
(\ref{drjj2}).

\vspace{5mm}

\noindent Next, we will give two $q$-analogues of
(\ref{recurmbes}):

\vspace{2mm}

{\bf Proposition 4.1.} {\sl The functions $J_{\nu}(x;q)$ and
${\cal J}_{\nu}(x;q)$
satisfy the (same) recurrence relations}
\begin{eqnarray}
\label{recj2}& &J_{\nu +m}(x;q) = R_{m,\nu}(x;q)J_{\nu}(x;q) -
R_{m-1,\nu +1}(x;q)J_{\nu -1}(x;q),\\[1ex]
\label{recjj2}& &{\cal J}_{\nu +m}(x;q) = R_{m,\nu}(x;q)
{\cal J}_{\nu}
(x;q) - R_{m-1,\nu +1}(x;q){\cal J}_{\nu -1}(x;q).
\end{eqnarray}

\vspace{2mm}

{\bf Proof:} We will give the proof for $J_{\nu}(x;q)$. We can
iterate
the recurrence relation (\ref{mix1}) to get an expression of the
form
\begin{eqnarray*}
J_{\nu + m}(x;q) = R_{m,\nu}(x;q)J_{\nu}(x;q) -
S_{m,\nu}(x;q)J_{\nu -1}(x;q).
\end{eqnarray*}
Next consider $J_{\nu + m + 1}(x;q)$. It satisfies the relations
\begin{eqnarray*}
& &J_{\nu +m+1}(x;q)= R_{m+1,\nu}(x;q)J_{\nu}(x;q) -
S_{m+1,\nu}(x;q)
J_{\nu-1}(x;q)\\[1ex]
& &\hspace{1cm}=R_{m,\nu+1}(x;q)J_{\nu+1}(x;q) - S_{m,\nu+1}(x;q)
J_{\nu}(x;q)\\[1ex]
& &\hspace{1cm}=\left[ R_{m,\nu+1}(x;q) \left\{ x +
\frac{1-q^{\nu}}{x}
\right\} -
S_{m,\nu+1}(x;q)\right] J_{\nu}(x;q) -
R_{m,\nu+1}(x;q)J_{\nu-1}(x;q).
\end{eqnarray*}
Hence $S_{m,\nu}(x;q) = R_{m-1,\nu+1}(x;q)$ and the proof is
completed.
Since $J_{\nu}(x;q)$ and
${\cal J}_{\nu}(x;q)$ satisfy the same recurrence relation (compare
(\ref{mix1}) and (\ref{recjj1})), it is obvious that (\ref{recjj2})
will also hold.$\Box$\newline

\noindent Another consequence of the proof is the recurrence
relation
\begin{eqnarray}
\label{recl1}R_{m+1,\nu}(x;q) = \left\{
x+\frac{1-q^{\nu}}{x}\right\}
R_{m,\nu+1}(x;q) - R_{m-1,\nu+2}(x;q).
\end{eqnarray}

\vspace{5mm}

Since (\ref{recj2}) and (\ref{recjj2}) are $q$-analogues of
(\ref{recurmbes}), the functions $R_{m,\nu}$ are $q$-analogues
of the Lommel polynomials. Later, when an explicit formula is
derived,
this will become clearer.

\vspace{5mm}

{\bf Lemma 4.2.} {\sl For non-integral $\nu$ the function
$R_{m,\nu}$
can be expressed as}
\begin{eqnarray}
\label{exlom1}R_{m,\nu}(x;q) &=& \frac{x e^{-i\pi\nu}
q^{-\frac{1}{2}
\nu(\nu-1)} (q;q)_{\infty}(q;q)_{\infty}}{(q^{\nu};q)_{\infty}
(q^{-\nu+1};q)_{\infty}}\times\\[1ex]
& &\hspace{2cm}\times\left\{ J_{\nu-1}(x;q){\cal J}_{\nu +m}(x;q)
-
J_{\nu +m}(x;q){\cal J}
_{\nu -1}(x;q) \right\}\nonumber.
\end{eqnarray}

\vspace{2mm}

{\bf Proof.} Multiply the recurrence relations (\ref{recj2}) and
(\ref
{recjj2}) with respectively ${\cal J}_{\nu-1}(x;q)$ and
$J_{\nu-1}(x;q)$
and subtract the resulting formulas to obtain
\begin{eqnarray*}
R_{m,\nu}(x;q) = \left\{ \frac{{\cal J}_{\nu+m}(x;q)J_{\nu-1}(x;q)
-
J_{\nu+m}(x;q){\cal J}_{\nu-1}(x;q)}{J_{\nu-1}(x;q){\cal
J}_{\nu}(x;q)
- J_{\nu}(x;q){\cal J}_{\nu -1}(x;q)} \right\}.
\end{eqnarray*}

Using (\ref{drj2}) and (\ref{drjj2}) we rewrite the denominator of
this
expression as
\begin{eqnarray*}
& &J_{\nu-1}(x;q){\cal J}_{\nu}(x;q) - J_{\nu}(x;q){\cal
J}_{\nu-1}(x;q)
\\[1ex]
&=& {\cal J}_{\nu}(x;q)\left[ \frac{1}{x}J_{\nu}(x;q) - \frac{1}{x}
q^{\frac{1}{2}\nu}J_{\nu}(xq^{\frac{1}{2}};q)\right] - J_{\nu}(x;q)
\left[ \frac{1}{x} {\cal J}_{\nu}(x;q) -
\frac{1}{x}q^{\frac{1}{2}\nu}
{\cal J}_{\nu}(xq^{\frac{1}{2}};q)\right]\\[1ex]
&=&\frac{q^{\frac{1}{2}\nu}}{x}\left\{ J_{\nu}(x;q){\cal J}_{\nu}(x
q^{\frac{1}{2}};q) - {\cal J}_{\nu}(x;q)J_{\nu}(xq^{\frac{1}{2}};q)
\right\}.
\end{eqnarray*}

The term in braces is a Wronskian type formula, which has been
evaluated by\newline
R.F. Swarttouw [9, (4.3.8)]. Explicitly,
\begin{eqnarray}
\label{wron}J_{\nu}(x;q){\cal J}_{\nu}(xq^{\frac{1}{2}};q) -
{\cal J}_{\nu}
(x;q)J_{\nu}(xq^{\frac{1}{2}};q) =
q^{\frac{1}{2}\nu(\nu-2)}e^{i\nu\pi}
\frac{(q^{\nu};q)_{\infty}(q^{-\nu+1};q)_{\infty}}{(q;q)_{\infty}
(q;q)_{\infty}},
\end{eqnarray}
which is non-zero for non-integral values of $\nu$. Now the lemma
follows easily.$\Box$

\vspace{2mm}

{\bf Remark:} The case $\nu=n\in\ZZ$ requires a closer
investigation,
since (\ref{wron}) equals zero if $\nu\in\ZZ$. However, for
$n\in\ZZ$
we have the relation, cf. [9, (3.2.12)]
\begin{eqnarray}
\label{minplus}J_n(x;q) = (-1)^n q^{-\frac{1}{2}n}
J_{-n}(xq^{-\frac{1}{2}n};q) = {\cal J}_n(x;q),
\end{eqnarray}
so that the numerator in (\ref{exlom1}) will also become  zero for
$n\in\ZZ$.
Now we can apply l'H\^{o}pital's rule to (\ref{exlom1}) to show
that the
formula still makes sense for $\nu\in\ZZ$. However, it will be
easier to
show that $R_{m,n}(x;q)$ is the analytic continuation of
$R_{m,\nu}(x;q)$,
when we have the explicit representation for the $R_{m,\nu}(x;q)$.
Until (\ref{explir}) we assume that $\nu$ is non-integral, but it
is
easily verified by analytic continuation that the results obtained
in
between are valid for integer $\nu$ as well.

\vspace{5mm}

Relation (\ref{exlom1}) gives rise to a three term recurrence
relation
for the $q$-Lommel polynomials:

\vspace{2mm}

{\bf Proposition 4.3.} {\sl The function $R_{m,\nu}(x;q)$ satisfies
the
recurrence relation}
\begin{eqnarray}
\label{recl2}\left\{ x^2 + 1 - q^{\nu+m} \right\} R_{m,\nu}(x;q) =
x \left\{ R_{m-1,\nu}(x;q) + R_{m+1,\nu}(x;q) \right\}.
\end{eqnarray}
with $R_{0,\nu}(x;q)=1$ and $R_{1,\nu}(x;q)=x+\frac{1-q^{\nu}}{x}$.

\vspace{2mm}

{\bf Proof:} $R_{0,\nu}(x;q)$ and $R_{1,\nu}(x;q)$ follow
immediately
from (\ref{recjj2}) and (\ref{recjj1}). In order to shorten the
notation
we introduce
\begin{eqnarray}
\label{cnu}C_{\nu}(x;q) =\frac{x
e^{-i\nu\pi}q^{-\frac{1}{2}\nu(\nu-1)}
(q;q)_{\infty}(q;q)_
{\infty}}{(q^{\nu};q)_{\infty}(q^{-\nu+1};q)_{\infty}}.
\end{eqnarray}
Now, starting with the lefthand side of (\ref{exlom1}), with $m$
replaced by $m+1$, we find with the recurrence relations
(\ref{mix1})
and (\ref{recjj1}), with $\nu$ replaced by $\nu+m$
\begin{eqnarray*}
R_{m+1,\nu}(x;q)&=&C_{\nu}(x;q) \left\{ J_{\nu-1}(x;q){\cal
J}_{\nu+m+1}
(x;q) - J_{\nu+m+1}(x;q){\cal J}_{\nu-1}(x;q)\right\}\\[1ex]
&=&C_{\nu}(x;q) \left\{J_{\nu-1}(x;q)\left[ \left\{ x +
\frac{1-q^{\nu
+m}}{x}\right\} {\cal J}_{\nu+m}(x;q) - {\cal
J}_{\nu+m-1}(x;q)\right]
\right.+ \\[1ex]
& &\hspace{1cm}-\left.{\cal J}_{\nu-1}(x;q)\left[
\left\{ x+\frac{1-q^{\nu+m}}
{x}\right\} J_{\nu+m}(x;q) - J_{\nu+m-1}(x;q) \right]
\right\}\\[1ex]
&=&\left\{ x + \frac{1-q^{\nu+m}}{x} \right\} R_{m,\nu}(x;q) -
R_{m-1,
\nu}(x;q)
\end{eqnarray*}
which proves the proposition.$\Box$

\vspace{5mm}

From Proposition 4.3 it is easy to see that the function
$R_{m,\nu}$
is not a polynomial, but a function of the form
\begin{eqnarray*}
R_{m,\nu}(x;q) = \sum\limits_{n=-m}^m c_n x^n.
\end{eqnarray*}
So $R_{m,\nu}$ is in fact a Laurent polynomial.
This suggests the introduction of the function
$p_m(x;q)=x^m R_{m,\nu}(x;q)$.
From (\ref{recl2}) it follows that $p_m(x;q)$ satisfies the
three term recurrence relation
\begin{eqnarray}
\label{prec} & &\left( x^2 + 1 - q^{\nu+m} \right) p_m(x;q) =
p_{m+1}(x;q) + x^2 p_{m-1}(x;q),\\[1ex]
& &\hspace{1cm} p_0(x;q)=1 \hspace{5mm} \mbox{ and } \hspace{5mm}
p_1(x;q)=x^2+\left(1-q^{\nu}\right).\nonumber
\end{eqnarray}
It follows from (\ref{prec}) that $p_{m}$ is actually a polynomial
of
degree $m$ in $x^2$.
The polynomial $p_m$ does not satisfy the conditions of Favard's
theorem, cf. [10, II.3.2], so that it is not orthogonal with
respect
to a positive
measure. However,\, $p_m(x;q)$ possesses nice properties, which are
$q$-analogues of the well-known properties of the Lommel
polynomials.

\vspace{5mm}

First we will derive a generating function for the $q$-Lommel
polynomials. Multiply (\ref{prec}) with $t^{m+1}$ and sum from
$m=1$ to
$\infty$. When we introduce the generating function
\begin{eqnarray*}
G(x,t) = \sum\limits_{m=0}^{\infty} t^m p_m(x;q),
\end{eqnarray*}
we find
\begin{eqnarray*}
& &x^2t \left( G(x,t) - 1 \right) + t \left( G(x,t) - 1 \right) -
q^{\nu}t \left( G(x,qt) - 1 \right) \\[1ex]
& &\hspace{1cm}= G(x,t) - 1 - t\left(x^2 + 1 - q^{\nu}\right) +
x^2t^2 G(x,t),
\end{eqnarray*}
which simplifies to
\begin{eqnarray}
G(x,t) = \frac{1 - q^{\nu}t G(x,qt)}{(1-t)(1-x^2t)}.\label{geng}
\end{eqnarray}
When we iterate relation (\ref{geng}) and make use of the fact that
$G(x,0)=1$, we have
\begin{eqnarray*}
G(x,t) = \sum\limits_{j=0}^{\infty} \frac{q^{j\nu}
q^{\frac{1}{2}j(j-1)}
(-t)^j}{(t;q)_{j+1} (x^2t;q)_{j+1}}.
\end{eqnarray*}
This expression is valid for $t\neq q^{-p}, \, x^2t\neq q^{-p}, \,
p\in\ZZ_+$. So we have the following generating function for the
function $R_{m,\nu}(x;q)$:
\begin{eqnarray}
\label{genr}\sum\limits_{m=0}^{\infty} (xt)^m R_{m,\nu}(x;q) &=&
\sum\limits_{j=0}^{\infty} \frac{q^{j\nu} q^{\frac{1}{2}j(j-1)}
(-t)^j}
{(t;q)_{j+1} (tx^2;q)_{j+1}}\\[1ex]
&=& \frac{1}{(1-t)(1-x^2t)} \qhyp{2}{2}{q,0}{qt, qtx^2}{tq^{\nu}}.
\nonumber
\end{eqnarray}

\vspace{5mm}

The generating function (\ref{genr}) gives rise to an explicit
expression of $R_{m,\nu}(x;q)$.
Use the $q$-binomial theorem [3, (1.3.2)] twice to obtain
\begin{eqnarray*}
& &\frac{1}{(t;q)_{j+1}} =
\frac{(tq^{j+1};q)_{\infty}}{(t;q)_{\infty}}
= \qhyp{1}{0}{q^{j+1}}{-}{t},\hspace{5mm} \left| t \right| <
1,\\[1ex]
& &\frac{1}{(tx^2;q)_{j+1}} = \qhyp{1}{0}{q^{j+1}}{-}{tx^2},
\hspace{5mm} \left| tx^2 \right| < 1.
\end{eqnarray*}
This yields
\begin{eqnarray*}
\sum\limits_{m=0}^{\infty} (xt)^m R_{m,\nu}(x;q) =
\sum\limits_{j,k,n=0}
^{\infty} \frac{q^{j\nu} q^{\frac{1}{2}j(j-1)} (-1)^j (q^{j+1};q)_k
(q^{j+1};q)_n}{(q;q)_k (q;q)_n} x^{2n} t^{j+k+n}.
\end{eqnarray*}
Equating powers of $t$ yields
\begin{eqnarray*}
x^m R_{m,\nu}(x;q) &=& \sum\limits_{n=0}^m \frac{x^{2n}}{(q;q)_n}
\sum\limits_{j=0}^{m-n} \frac{q^{j\nu} q^{\frac{1}{2}j(j-1)} (-1)^j
(q^{j+1};q)_{m-n-j} (q^{j+1};q)_n}{(q;q)_{m-n-j}}\\[1ex]
&=& \sum\limits_{n=0}^m x^{2n} \sum\limits_{j=0}^{m-n}
\frac{(q^{n-m};q)_j (q^{n+1};q)_j}{(q;q)_j (q;q)_j} q^{j\nu}
q^{j(m-n)}
\\[1ex]
&=& \sum\limits_{n=0}^m x^{2n} \qhyp{2}{1}{q^{n-m},q^{n+1}}{q}
{q^{\nu+m-n}}.
\end{eqnarray*}
This expression for $R_{m,\nu}(x;q)$ shows the analytic dependence
on
$q^{\nu}$, cf. the remark following Lemma 4.2.
Finally, using Heine's transformation formula, cf. [3, (1.4.1)], we
obtain the explicit representation
\begin{eqnarray}
\label{explir}R_{m,\nu}(x;q) = \sum\limits_{n=0}^m x^{2n-m} \frac
{(q^{n+1};q)_{\infty}(q^{\nu};q)_{\infty}}{(q;q)_{\infty}
(q^{\nu+m-n};q)_{\infty}} \qhyp{2}{1}{q^{-n},q^{\nu+m-n}}{q^{\nu}}
{q^{n+1}}.
\end{eqnarray}
The explicit expression (\ref{explir}) can also be obtained from
(\ref{exlom1}) in combination with the explicit series
representations
(\ref{calj}) and (\ref{heq}). In this case a formula [9, (6.4.4)]
for
the product of two Hahn-Exton $q$-Bessel functions has to be used.

\vspace{5mm}

Formally we can obtain a $q$-analogue of Hurwitz's formula
(\ref{hur})
by taking termwise limits in the explicit representation
(\ref{explir}).
This gives
\begin{eqnarray*}
& &x^m R_{m,\nu}(x;q)
\stackrel{m\rightarrow\infty}{\longrightarrow}
\frac{(q^{\nu};q)_{\infty}}{(q;q)_{\infty}}
\sum\limits_{n=0}^{\infty}
x^{2n} (q^{n+1};q)_{\infty} \qhyp{2}{1}{q^{-n},0}{q^{\nu}}{q^{n+1}}
\\[1ex]
& &\hspace{3cm}= (q^{\nu+1};q)_{\infty} \sum\limits_{k=0}^{\infty}
\frac{q^k}{(q^{\nu};q)_k
(q;q)_k} \sum\limits_{n=k}^{\infty} \frac{(q^{-n};q)_k}{(q;q)_n}
q^{nk} x^{2n}.
\end{eqnarray*}
The inner sum equals
\begin{eqnarray*}
(-1)^k q^{\frac{1}{2}k(k-1)} x^{2k}
\frac{1}{(x^2;q)_{\infty}},
\end{eqnarray*}
so that we formally obtain
\begin{eqnarray}
x^m R_{m,\nu}(x;q) &\stackrel{m\rightarrow\infty}{\longrightarrow}&
\frac{(q^{\nu};q)_{\infty}}{(x^2;q)_{\infty}}
\qhyp{1}{1}{0}{q^{\nu}}
{qx^2}\nonumber\\[1ex]
&=& \frac{(q;q)_{\infty}}{(x^2;q)_{\infty}} x^{1-\nu}
J_{\nu-1}(x;q).
\end{eqnarray}

\vspace{5mm}

\subsection{$q$-Lommel polynomials associated with (2.14)}

In this subsection we consider the $q$-analogues of the Lommel
polynomials that arise from the iteration of (\ref{mix2}). Note
that
the shift in the argument in (\ref{mix2}) will present some
difficulty. As a consequence, there is no unique definition of
these
$q$-Lommel polynomials. For a certain choice we show that the
(modified)
$q$-Lommel polynomials associated with the Jackson $q$-Bessel
function, see M.E.H. Ismail \cite{Is1}, are obtained in this
context.
These polynomials are orthogonal as shown by M.E.H. Ismail [4,
(4.16)].
For a slightly different definition of the $q$-Lommel polynomials
we
present a Hurwitz-type formula, cf. (\ref{hur}).

We start with iterating (\ref{mix2}).

\vspace{5mm}

{\bf Lemma 4.4.} {\sl There exist unique constants $a_i(\nu,m)$ and
$b_j(\nu,m)$ so that
\begin{eqnarray*}
J_{\nu+m}(xq^m;q^2) = \sum\limits_{i=0}^{[m/2]} \frac{a_i(\nu,
m)}{x^{m-2i}} J_{\nu}(xq^i;q^2) + \sum\limits_{j=0}^
{[(m-1)/2]} \frac{b_j(\nu,m)}{x^{m-1-2j}} J_{\nu-1}(xq^j;q^2),
\end{eqnarray*}
where $[a]$ denotes the largest integer smaller than or equal to
$a$.}

\vspace{2mm}

{\bf Proof:} The proof uses induction with respect to $m$. The case
$m=0$ is trivial and the case $m=1$ is just (\ref{mix2}). For the
induction step we use, cf. (\ref{mix2}),
\begin{eqnarray}
& &J_{\nu+m+1}(xq^{m+1};q^2) = J_{(\nu+m)+1}(xq^m\cdot q;q^2)
\label{ind}\\[1ex]
& &\hspace{1cm}= q^{-(\nu+m+1)} \frac{1-q^{2(\nu+m)}}{xq^m}
J_{\nu+m}(xq^m;q^2) -
q^{-(\nu+m+1)} J_{\nu+(m-1)}(xq\cdot q^{m-1};q^2).\nonumber
\end{eqnarray}
Use the induction hypothesis in the last two terms of (\ref{ind})
and
shift summation parameters to prove the lemma. $\Box$

\vspace{2mm}

Actually this proof yields a recurrence relation for the constants
$a_i(\nu, m+1)$. Explicitly,
\begin{eqnarray}
\label{reca}a_i(\nu,m+1) =
q^{-(\nu+1)}q^{-2m}(1-q^{2(\nu+m)})a_i(\nu,m)
- q^{-\nu + 2(i-1-m)} a_{i-1}(\nu,m-1),
\end{eqnarray}
with $a_i(\nu,m)=0$ for $i < 0$ or $i > [m/2]$.
A similar result holds for the constants $b_j(\nu, m)$, but these
constants will be related to the $a_i(\nu, m)$ in due course, cf.
Lemma
4.5. So we do not present the recurrence relation for the
$b_j(\nu, m)$'s.

\vspace{5mm}

Although Lemma 4.4 also holds with $J_{\nu}$ replaced by
${\cal J}_{\nu}$, cf. (\ref{calj}), this does not lead to nice
expressions as in the previous subsection due to the shift in the
argument in Lemma 4.4.

\vspace{5mm}

Motivated by Lemma 4.4 and the recurrence relation (\ref{reca}) we
define
\begin{eqnarray}
r_{m,\nu}(x;q^2) = \sum\limits_{i=0}^{[m/2]} \frac{a_i(\nu,m)}
{x^{m-2i}} q^{-i(i+1)},\label{rklein}
\end{eqnarray}
which is a polynomial of degree $m$ in $\frac{1}{x}$. The factor
$q^{-i(i+1)}$ in (\ref{rklein}) is chosen in such a way that
(\ref{reca}) yields a three term recurrence relation for the
$r_{m,\nu}(x;q^2)$. Explicitly
\begin{eqnarray}
\label{recr2}q^{2m} r_{m+1,\nu}(x;q^2) = q^{-(\nu+1)}
\frac{1-q^{2(\nu+m)}}{x} r_{m,\nu}(x;q^2) - q^{-(\nu+2)}
r_{m-1,\nu}(x;q^2),
\end{eqnarray}
which follows from multiplying (\ref{reca}) by
$x^{-m-1+2i} q^{-i(i+1)}$
and summing from $i=0$ to $[(m+1)/2]$. The begin values of the
recurrence relation (\ref{recr2}) are
\begin{eqnarray}
\label{begrec}r_{0,\nu}(x;q^2) = 1, \hspace{1cm} r_{1,\nu}(x;q^2)
=
q^{-(\nu+1)} \left( 1-q^{2\nu} \right) x^{-1}.
\end{eqnarray}

\vspace{2mm}

In order to determine the constants $a_i(\nu, m)$ in (\ref{rklein})
we introduce the polynomial $h_{m,\nu}(x;q^2) = r_{m,\nu}\left(
\frac{1}{x};q^2\right)$. The three term recurrence relation
(\ref{reca})
with its begin values (\ref{begrec}) is rewritten as
\begin{eqnarray}
& &q^{2m} h_{m+1,\nu}(x;q^2) = q^{-\nu-1} \left( 1-q^{2(\nu+m)}
\right)
x h_{m,\nu}(x;q^2) - q^{-\nu-2} h_{m-1,\nu}(x;q^2),\nonumber\\[1ex]
& &\hspace{1cm} h_{0,\nu}(x;q^2) = 1 \hspace{5mm} \mbox{ and }
\hspace{5mm} h_{1,\nu}(x;q^2) = q^{-\nu-1} \left(1-q^{2\nu} \right)
x.
\label{modrec}
\end{eqnarray}

\vspace{5mm}

Let us recall the (modified) $q$-Lommel polynomials associated with
the
Jackson $q$-Bessel function as presented by M.E.H. Ismail
\cite{Is1}.
Define, cf. [4, (3.6)],
\begin{eqnarray}
\tilde{h}_{m,\nu}(x;q^2) = \sum\limits_{j=0}^{[m/2]}
\frac{(2x)^{m-2j}
(-1)^j (q^{2\nu};q^2)_{m-j} (q^2;q^2)_{m-j}}{(q^2;q^2)_j
(q^{2\nu};q^2)_j (q^2;q^2)_{m-2j}}q^{2j(j+\nu-1)},\label{his}
\end{eqnarray}
then the set of polynomials $\tilde{h}_{m,\nu}(x;q^2)$ satisfies
the
three term recurrence relation, cf. [4, (1.22)],
\begin{eqnarray}
& &\tilde{h}_{m+1,\nu}(x;q^2) = 2x \left( 1-q^{2(\nu+m)} \right)
\tilde{h}_{m,\nu}(x;q^2) - q^{2(m+\nu-1)}
\tilde{h}_{m-1,\nu}(x;q^2),
\nonumber\\[1ex]
& &\hspace{1cm} \tilde{h}_{0,\nu}(x;q^2) = 1 \hspace{5mm} \mbox{
and }
\hspace{5mm} \tilde{h}_{1,\nu}(x;q^2) = 2x \left(1-q^{2\nu}
\right).
\label{modis}
\end{eqnarray}
A straightforward calculation shows that (\ref{modrec}) and
(\ref{modis}) are equivalent. It suffices to take
\begin{eqnarray}
h_{m,\nu}(x;q^2) = q^{-\frac{3}{2}m\nu - m^2} \tilde{h}_{m,\nu}
(q^{\frac{\nu}{2}}x/2;q^2).\label{hish}
\end{eqnarray}
From (\ref{hish}), (\ref{his}) and (\ref{rklein}) we obtain the
explicit
expression
\begin{eqnarray}
a_i(m,\nu) = q^{-m(m+\nu)} q^{i(3i+\nu-1)} \frac{(-1)^i
(q^{2\nu};q^2)_{m-i} (q^2;q^2)_{m-i}}{(q^2;q^2)_i (q^{2\nu};q^2)_i
(q^2;q^2)_{m-2i}}.\label{expla}
\end{eqnarray}

\vspace{5mm}

Now that we have determined the value of $a_i(\nu,m)$ we present
its
link to the constants $b_j(\nu,m)$.

\vspace{2mm}

{\bf Lemma 4.5.}
\begin{eqnarray}
b_j(\nu,m+1) = -q^{2j-m-\nu-1} a_j(\nu+1,m).
\end{eqnarray}

\vspace{2mm}

{\bf Proof:} Consider
\begin{eqnarray}
J_{\nu+m+1}(xq^{m+1};q^2) = J_{(\nu+1)+m}(xq
\cdot q^m;q^2)\label{lem45}
\end{eqnarray}
and apply Lemma 4.4 with $\nu$ replaced by $\nu+1$ to the right
hand
side of (\ref{lem45}). In the resulting sum involving the
Hahn-Exton
$q$-Bessel function of order $\nu+1$ we use (\ref{mix2}).
Comparison
of the coefficients of $x^{2i-m} J_{\nu-1}(xq^i;q^2)$ with the ones
obtained by applying Lemma 4.4 with $m$ replaced by $m+1$ to the
left
hand side of (\ref{lem45}) proves the lemma.$\Box$

\vspace{2mm}

In the proof of Lemma 4.5 we can also compare the coefficients of
$x^{2i-m-1} J_{\nu}(xq^i;q^2)$ and this results in
\begin{eqnarray}
\label{abcom}a_i(\nu,m+1) = q^{i-m-\nu-1} \left( 1-q^{2\nu} \right)
a_i(\nu+1,m) + q^{2i-m-1} b_{i-1}(\nu+1,m).
\end{eqnarray}
Using Lemma 4.5 and (\ref{rklein}) we can rewrite (\ref{abcom}) as
\begin{eqnarray}
\label{recrecr}r_{m+1,\nu}(x;q^2) &=&
\frac{1}{x}\left(1-q^{2\nu}\right)
q^{-(\frac{1}{2}m+\nu+1)} r_{m,\nu+1}(xq^{\frac{1}{2}};q^2)\\[1ex]
& &\hspace{3cm}- q^{-(m+\nu+3)} r_{m-1,\nu+2}(xq;q^2).\nonumber
\end{eqnarray}
Since we have identified $r_{m,\nu}(x;q^2)$ with $q$-analogues of
Lommel polynomials associated with Jackson's $q$-Bessel functions,
Hurwitz's formula, cf. M.E.H. Ismail [4, (3.10)], holds for
$r_{m,\nu}(x;q^2)$. We want to have a $q$-analogue of Hurwitz's
formula involving the Hahn-Exton $q$-Bessel function and in order
to
do so we make use of the following $q$-analogue of the Lommel
polynomial:
\begin{eqnarray}
\tilde{r}_{m,\nu}(x;q) = \sum\limits_{i=0}^{[m/2]}
\frac{a_i(m,\nu)}
{x^{m-2i}} q^{-2i(i-1)}.\label{tilder}
\end{eqnarray}

From (\ref{reca}) we see that it satisfies the recurrence relation
\begin{eqnarray}
\tilde{r}_{m+1,\nu}(x;q^2) = x^{-1} q^{-\nu-1-2m}
\left(1-q^{2(\nu+m)}
\right) \tilde{r}_{m,\nu}(x;q^2) - q^{1-\nu-3m} \tilde{r}_{m-1,\nu}
(q^{-1}x;q^2).\label{rectilder}
\end{eqnarray}
It follows from Favard's theorem, cf. [10, II.3.2], that
$\tilde{r}_{m,\nu}\left(\frac{1}{x};q^2\right)$ does not yield a
set
of orthogonal polynomials. However, from (\ref{expla}) and
(\ref{tilder}) we get
\begin{eqnarray}
\tilde{r}_{m,\nu}(x;q^2) &=& x^{-m}q^{-m(m+\nu)} (q^{2\nu};q^2)_m
\label{exrtilde}\\[1ex]
&\times&\qqhyp{5}{3}
{q^{1-m}, -q^{1-m}, q^{-m},
-q^{-m}, 0}{q^{2\nu}, q^{2(1-m-\nu)}, q^{-2m}}
{x^2q^{2-\nu}}\nonumber
\end{eqnarray}
and from (\ref{exrtilde}) we obtain the following $q$-analogue of
Hurwitz's formula
\begin{eqnarray*}
x^m q^{m(m+\nu)} \tilde{r}_{m,\nu}(x;q^2)
\stackrel{m\rightarrow\infty}
{\longrightarrow}
(q^2;q^2)_{\infty} x^{1-\nu} q^{\frac{1}{2}\nu(1-\nu)}
J_{\nu-1}(xq^{\frac{\nu}{2}};q^2),
\end{eqnarray*}
by taking formal termwise limits in (\ref{exrtilde}).

\end{document}